\documentclass[12pt,a4paper,draft]{article}
\usepackage[english,russian]{babel}
\usepackage{amssymb, amsmath, latexsym, amsfonts, amsthm}
\begin{document}
\renewcommand{\refname}{References}
\newtheorem{lemma}{Lemma}
\newtheorem*{theorem}{Theorem}
\newtheorem*{theorem*}{Theorem A}
\newtheorem*{theorem**}{Theorem B}
\newtheorem{definition}{Definition}
\newtheorem{corollary}{Corollary}
\newtheorem*{corollary*}{Corollary A}
\newtheorem{proposition}{Proposition}
\newtheorem*{proposition M}{Proposition M}
\begin{center}
ON A HILBERT SPACE OF ENTIRE FUNCTIONS
\end{center}
\begin{center}
Il'dar Kh. Musin\footnote {The research was supported by RFBR (15-01-01661) and Program of the Presidium of RAS (project "`Complex Analysis and Functional Equations"').}
\end{center}
\begin{center}
Institute of Mathematics with Computer Centre of Ufa Scientific Centre of Russian Academy of Sciences, 
Chernyshevsky str., 112, Ufa, 450077, Russia
\end{center}

\vspace {0.1 cm}

\renewcommand{\abstractname}{}
\begin{abstract}
{\bf Abstract:} 
A weighted Hilbert space $F^2_{\varphi}$ of entire functions of $n$ variables is considered in the paper. The weight function $\varphi$ is a convex function on ${\mathbb C}^n$ depending on modules of variables and growing at infinity faster than  $a \Vert z \Vert$ for each $a > 0$. The problem of description of the strong dual of this space in terms of the Laplace transformation of functionals is studied in the article. Under some additional conditions on $\varphi$ the space of the Laplace transforms of linear continuous functionals on $F^2_{\varphi}$ is described. The proof of the main result is based on new properties of the Young-Fenchel transformation and a result of R.A. Bashmakov, K.P. Isaev and R.S. Yulmukhametov on asymptotics of multidimensional Laplace transform.

\vspace {0.1 cm}


\vspace {0.1 cm}


\end{abstract}

\vspace {0.1 cm}

\begin{section}
{\bf  Introduction} 
\end{section}

{\bf 1.1. On the Problem}. 
Let $H({\mathbb C}^n)$ be the space of entire functions in ${\mathbb C}^n$, $d \mu_n$ be the Lebesgue measure in ${\mathbb C}^n$ and for $u=(u_1, \ldots , u_n) \in {\mathbb R}^n \ ({\mathbb C}^n)$ \ $abs \, u: = (\vert u_1 \vert, \ldots , \vert u_n \vert)$. 

Denote by ${\mathcal V}({\mathbb R}^n)$ the set of all convex functions $g$ in 
${\mathbb R}^n$ such that:

1). $g(x_1, \ldots , x_n) = g(\vert x_1 \vert, \ldots , \vert x_n \vert), \ (x_1, \ldots , x_n)\in {\mathbb R}^n$;

2). the restriction $g$ on $[0, \infty)^n$ is nondecreasing in each variable;

3). $\displaystyle \lim_{x \to \infty} \frac {g(x)}{\Vert x \Vert}= + \infty$ \ ($\Vert x \Vert$ is the Euclidean norm of $x \in {\mathbb R}^n$).  

To each $\varphi \in {\mathcal V}({\mathbb R}^n)$ we associate the Hilbert space
$$
F^2_{\varphi} = \{f \in H({\mathbb C}^n): \Vert f \Vert_{\varphi} = 
\left(\int_{{\mathbb C}^n} \vert f(z)\vert^2 
e^{- 2 \varphi (abs \ z)} \ d \mu_n (z)\right)^{\frac 1 2} < \infty \}
$$
with the scalar product
$$
(f, g)_{\varphi} = \int_{{\mathbb C}^n} f(z) \overline {g(z)} e^{- 2 \varphi (abs \ z)} \ d \mu_n (z), \ f, g \in F^2_{\varphi}.
$$
If $\varphi (x) = \frac {\Vert x \Vert^2}{2}$ then $F^2_{\varphi}$ is the Fock space.

Obviously, for each function $\varphi \in {\mathcal V}({\mathbb R}^n)$  and each $\lambda \in {\mathbb C}^n$ the function 
$f_{\lambda}(z) = e^{\langle \lambda, z \rangle}$ is in $F^2_{\varphi}$. 
So for each linear continuous functional $S$ on $F^2_{\varphi}$ the function
$$
\hat S (\lambda) = S(e^{\langle \lambda, z \rangle}), \ \lambda \in {\mathbb C}^n, 
$$
is well defined in ${\mathbb C}^n$. $\hat S$ is called  the Laplace transform of the functional $S$.  It is easy to see that $\hat S$ is an entire function.

Denote the dual space for $F^2_{\varphi}$ by $(F^2_{\varphi})^*$. 

The aim of the work is to find the conditions for $\varphi \in {\mathcal V}({\mathbb R}^n)$ allowing to describe the space $\widehat {(F^2_{\varphi})^*}$ of
the Laplace transforms of the linear continuous functionals on $F^2_{\varphi}$ as $F^2_{\varphi^*}$.

If $\varphi (x) = \frac {\Vert x \Vert^2}{2}$ then
$\widehat {(F^2_{\varphi})^*} = F^2_{\varphi}$. Indeed, in this case the problem of describing of the space $(F^2_{\varphi})^*$ in terms of the Laplace transform of the functionals is easily solved thanks to the classical
representation: for each $f \in F^2_{\varphi}$
$$
f(\lambda) = \pi^{-n} \int_{{\mathbb C}^n} f(z) 
e^{\langle \lambda, \overline z \rangle - \Vert z \Vert^2} \ d \mu_n (z), \ \lambda \in {\mathbb C}^n.
$$

In case of a radial function $\varphi \in {\mathcal V}({\mathbb R}^n)$ the above  mentioned problem was solved by V.V. Napalkov and S.V. Popenov \cite {N-P}, \cite {P}.  

{\bf 1.2. Notations and definitions}. For $u=(u_1, \ldots , u_n)$, $v=(v_1, \ldots , v_n) \in {\mathbb R}^n ({\mathbb C}^n)$ let $\langle u, v \rangle := u_1 v_1 + \cdots + u_n v_n$,
$\Vert u \Vert$ be the Euclidean norm of $u$.

For $\alpha = (\alpha_1, \ldots , \alpha_n) \in {\mathbb Z}_+^n$, 
$z =(z_1, \ldots , z_n) \in {\mathbb C}^n$ let $\vert \alpha \vert := \alpha_1 + \ldots  + \alpha_n$, $\tilde \alpha: = (\alpha_1 + 1, \ldots , \alpha_n + 1)$,
$z^{\alpha} := z_1^{\alpha_1} \cdots z_n^{\alpha_n}$, 
$D^{\alpha}_z :=
\frac {{\partial}^{\vert \alpha \vert}}{\partial z_1^{\alpha_1} \cdots \partial z_n^{\alpha_n}}$.

For $\alpha = (\alpha_1, \ldots , \alpha_n) \in {\mathbb Z}_+^n$, 
$\varphi \in {\mathcal V}({\mathbb R}^n)$  let
$$
c_{\alpha} (\varphi):= 
\int_{{\mathbb C}^n} \vert z_1 \vert^{2 \alpha_1} \cdots \vert z_n \vert^{2 \alpha_n} 
e^{- 2 \varphi (abs \, z)} \, d \mu_n (z).
$$

For a function $u$ with a domain containing the set $(0, \infty)^n$ define a function
$u[e]$ on ${\mathbb R}^n$ by the rule:
$$
u[e](x) = u(e^{x_1}, \ldots, e^{x_n}), \ x = (x_1, \ldots , x_n) \in  {\mathbb R}^n.
$$
  
By ${\mathcal B}({\mathbb R}^n)$ denote the set of all continuous functions 
$u: {\mathbb R}^n \to {\mathbb R}$ satisfying the condition 
$\displaystyle \lim_{x \to \infty} \frac {u(x)}{\Vert x \Vert}=~+\infty$. 

The Young-Fenchel conjugate of a function $u:{\mathbb R}^n \to [-\infty, + \infty]$ is the function 
$u^*:{\mathbb R}^n \to [-\infty, + \infty]$ defined by
$$
u^*(x) = \displaystyle \sup \limits_{y \in {\mathbb R}^n}(\langle x, y \rangle - u(y)), \ x \in {\mathbb R}^n. 
$$

If $E$ is a convex domain in ${\mathbb R}^n$, $h$ is a convex function on $E$, $\tilde E =\{y \in {\mathbb R}^n: h^*(y) < \infty \}$, $p > 0$, then 
$
D_y^h (p): = \{x \in E: h(x) + h^*(y) - \langle x, y \rangle \le p \}, \ y \in \tilde E.
$

Denote the $n$-dimensional volume of a set $D \subset {\mathbb R}^n$ by $V(D)$.

{\bf 1.3. Main result}. 

\begin{theorem} 
Let $\varphi \in {\mathcal V}({\mathbb R}^n)$ 
and for some $K > 0$ 
$$
\frac 1 K \le V(D_{\alpha}^{\varphi [e]}(1 / 2)) V(D_{\alpha}^{\varphi^* [e]}(1 / 2)) 
\prod \limits_{j =1}^n \alpha_j \le K, \ \forall \alpha  = (\alpha_1, \ldots,  \alpha_n) \in {\mathbb N}^n.
$$

Then the mapping ${\mathcal L}: S \in (F^2_{\varphi})^* \to \hat S$ establishes an isomorphism between the spaces $(F^2_{\varphi})^*$ and $F^2_{\varphi^*}$. 
\end{theorem}

The proof of the Theorem in Subsection 3.2 is based on new properties of the Young-Fenchel transform (Subsection 2.1) and one result on the asymptotics of the multidimensional Laplace integral in \cite {BIY} (Subsection 2.2). 

\section{Auxiliary results}

2.1. {\bf On some properties of the Young-Fenchel transform}. It is easy to verify that the following statement holds.

\begin{proposition}
Let $u \in {\mathcal B}({\mathbb R}^n)$. Then $(u[e])^*(x) > -\infty$ for $x \in {\mathbb R}^n$, $(u[e])^*(x)=+\infty$ for $x \notin [0, \infty)^n$, $(u[e])^*(x) <+\infty$ for
$x \in [0, \infty)^n$.
\end{proposition}

We only note that the last statement of Proposition 1 is implied, for instance, by the fact that for each $M>0$ there exists a constant $A>0$ such that 
$$
(u[e])^*(x) \le 
\sum \limits_{1 \le j \le n: x_j \ne 0} (x_j \ln\frac {x_j}{M} - x_j) + A, \ x \in [0, \infty)^n.
$$ 

\begin{proposition}
Let $u \in {\mathcal B}({\mathbb R}^n)$. Then 
$$
\displaystyle \lim_{x \to \infty, \atop x \in [0, \infty)^n} \frac {(u[e])^*(x)}{\Vert x \Vert}= + \infty. 
$$
\end{proposition}

{\bf Proof}. For any $x \in [0, \infty)^n$ and $t \in {\mathbb R}^n$ we have that
$$
(u[e])^*(x) \ge \langle x, t \rangle - (u[e])(t). 
$$
Using this inequality we obtain that for each $M >0$ 
$$
(u[e])^*(x) \ge M \Vert x \Vert - u[e]\left(\frac {M x}{\Vert x \Vert}\right), \ x \in [0, \infty)^n \setminus \{0\}.
$$ 
From this the assertion follows. $\square$

The next three statements were proved in \cite {M} (Lemma 6, Proposition 3, Proposition 4).

\begin{proposition}
Let $u \in {\mathcal B}({\mathbb R}^n)$.
Then
$$
(u[e])^*(x) + (u^*[e])^*(x) \le \sum \limits _{1 \le j \le n: \atop x_j \ne 0}
(x_j \ln x_j - x_j), \
x = (x_1, \ldots , x_n) \in [0, \infty)^n \setminus \{0\};
$$
$$
(u[e])^*(0) + (u^*[e])^*(0) \le 0.
$$
\end{proposition}

\begin{proposition}
Let $u \in {\mathcal B}({\mathbb R}^n) \cap C^2({\mathbb R}^n)$ be a convex function. Then
$$
(u[e])^*(x) + (u^*[e])^*(x) = \sum \limits_{j =1}^n
(x_j \ln x_j - x_j), \
x = (x_1, \ldots , x_n)\in (0, \infty)^n.
$$
\end{proposition}

\begin{proposition}
Let $u \in {\mathcal V}({\mathbb R}^n) \cap C^2({\mathbb R}^n)$. Then
$$
(u[e])^*(x) + (u^*[e])^*(x) =  \sum \limits _{1 \le j \le n: \atop x_j \ne 0} (x_j \ln x_j - x_j), \,
x = (x_1, \ldots , x_n) \in [0, \infty)^n \setminus \{0\};
$$
$$
(u[e])^*(0) + (u^*[e])^*(0) = 0.
$$
\end{proposition}

Propositions 4 and 5 can be strengthen using the results by D. Azagra \cite {A 1}, \cite {A 2}. He proved the following theorem.

\begin{theorem*} 
Let $U \subseteq {\mathbb R}^n$ be an open convex set. For each convex function $f: U \to  {\mathbb R}$ and each $\varepsilon > 0$ there exists a real analytic convex function $g: U \to  {\mathbb R}$ such that 
$$
f(x) - \varepsilon \le g(x) \le f(x), \ x \in U.
$$
\end{theorem*} 

Thus, the following corollary holds \cite {A 2}.

\begin{corollary*}
Let $U \subseteq {\mathbb R}^n$ be an open convex set. 
For each convex function $f: U \to  {\mathbb R}$ and each $\varepsilon > 0$ there exists  an infinitely differentiable convex function $g: U \to  {\mathbb R}$ such that 
$f(x) - \varepsilon \le g(x) \le f(x), \ x \in U$.
\end{corollary*}

Using Proposition 4 and Corollary A we easily confirm the following statement.

\begin{proposition}
Let $u \in {\mathcal B}({\mathbb R}^n)$ be a convex function. Then
$$
(u[e])^*(x) + (u^*[e])^*(x) = \sum \limits_{j =1}^n 
(x_j \ln x_j - x_j), \
x = (x_1, \ldots , x_n)\in (0, \infty)^n.
$$
\end{proposition}

Moreover, the following proposition is true.

\begin{proposition}
Let $u \in {\mathcal V}({\mathbb R}^n)$ be a convex function. Then
$$
(u[e])^*(x) + (u^*[e])^*(x) =  \sum \limits _{1 \le j \le n: \atop x_j \ne 0} (x_j \ln x_j - x_j), \,
x = (x_1, \ldots , x_n) \in [0, \infty)^n \setminus \{0\};
$$
$$
(u[e])^*(0) + (u^*[e])^*(0) = 0.
$$
\end{proposition}

{\bf Proof}. By Proposition 6 our statement is true for $x \in (0, \infty)^n$. 
Now let $x =(x_1, \ldots , x_n)$ belongs to the boundary of $[0, \infty)^n$ and $x \ne 0$. 
For the sake of simplicity we consider the case when the first $k$ ($1 \le k \le n-1$) coordinates of $x$ are positive and all other are equal to zero. 
For all $\xi = (\xi_1, \ldots , \xi_n), \mu = (\mu_1, \ldots , \mu_n) \in {\mathbb R}^n$ we have that
$$
(u[e])^*(x) + (u^*[e])^*(x) 
\ge \sum \limits _{j=1}^k x_j (\xi_j + \mu_j) - 
(u(e^{\xi_1}, \ldots , e^{\xi_n}) + u^*(e^{\mu_1}, \ldots , e^{\mu_n})).
$$
From this inequality we obtain that
$$
(u[e])^*(x) + (u^*[e])^*(x) 
\ge \sum \limits _{j=1}^k x_j (\xi_j + \mu_j) - 
$$
$$
-
(u(e^{\xi_1}, \ldots , e^{\xi_k}, 0, \ldots , 0) + 
u^*(e^{\mu_1}, \ldots , e^{\mu_k}, 0, \ldots , 0)).
$$
Define a function
$u_k$ on ${\mathbb R}^k$ by the rule: 
$(\lambda_1, \ldots , \lambda_k) \in {\mathbb R}^k \to 
u(\lambda_1, \ldots , \lambda_k, 0, \ldots , 0)$.
Note that for each $t= (t_1, \ldots , t_k) \in {\mathbb R}^k, 
\breve t = (t_1, \ldots , t_k, 0, \ldots , 0)  \in {\mathbb R}^n$
$$
u^*(\breve t) = 
\sup \limits_{v \in {\mathbb R}^n} (\langle \breve t, v \rangle - u(v)) \le 
$$
$$
\le
\sup \limits_{v_1, \ldots , v_k \in {\mathbb R}} 
(\sum \limits_{j=1}^k t_j v_j   - u(v_1, \ldots , v_k, 0, \ldots , 0)) = 
\sup \limits_{v \in {\mathbb R}^k} (\langle t, v \rangle - u_k(v)) = u_k^*(t).
$$
Using this and the above inequality we have that for $\tilde x = (x_1, \ldots , x_k) \in {\mathbb R}^k$ and all
$\tilde \xi = (\xi_1, \ldots , \xi_k), 
\tilde \mu = (\mu_1, \ldots , \mu_k) \in {\mathbb R}^k$ 
$$
(u[e])^*(x) + (u^*[e])^*(x) 
\ge \langle \tilde x,  \tilde \xi \rangle - u_k[e](\tilde \xi)  + 
\langle \tilde x,  \tilde \mu \rangle - u_k^*[e](\tilde \mu).
$$
Hence,
$$
(u[e])^*(x) + (u^*[e])^*(x) 
\ge (u_k[e])^*(\tilde x) + (u_k^*[e])^*(\tilde x). 
$$
Since by the Proposition 6  
$$
(u_k[e])^*(\tilde x) + (u_k^*[e])^*(\tilde x) = 
\sum \limits_{j =1}^k(x_j \ln x_j - x_j),
$$
then 
$
(u[e])^*(x) + (u^*[e])^*(x) \ge 
\sum \limits_{j =1}^k(x_j \ln x_j - x_j). 
$
From this and Proposition 3 the first statement of the Proposition follows.

If $x =0$ then
$(u[e])^*(0) = - \inf \limits_{\xi \in {\mathbb R}^n} u[e] (\xi) = - u(0)$, 
$(u^*[e])^*(0) = - \inf \limits_{\xi \in {\mathbb R}^n} u^*[e] (\xi) = -u^*(0) = \inf \limits_{\xi \in {\mathbb R}^n} u(\xi) = u(0)$. 
Therefore, 
$(u[e])^*(0) + (u^*[e])^*(0) = 0$. $\square$

2.2. {\bf Asymptotics of multidimensional Laplace integral}. In \cite {BIY} there was established the following theorem.

\begin{theorem**} 
Let $E$ be a convex domain in ${\mathbb R}^n$, $h$ be a convex function on $E$, $\tilde E =\{y \in {\mathbb R}^n: h^*(y) < \infty \}$ and the interior of $\tilde E$ is not empty. Let
$$
D^h = \{(x, y) \in {\mathbb R}^n \times {\mathbb R}^n: 
h(x) + h^*(y) - \langle x, y \rangle \le 1 \},
$$
$$
D_y^h = \{x \in {\mathbb R}^n: (x, y) \in D \}, \ y \in {\mathbb R}^n.
$$
Then
$$
e^{-1} V(D_y^h) e^{h^*(y)} \le 
\int_{{\mathbb R}^n} e^{\langle x, y \rangle - h(x)} \, d x  
\le (1 + n!) V(D_y^h) e^{h^*(y)}, \ y \in \tilde E.
$$
\end{theorem**} 

Here it is assumed that $h(x) = +\infty$ for $x \notin E$.

\section{Description of the dual space}

{\bf 3.1. Auxiliary lemmas}. In the proof of the Theorem the following four lemmas will be useful.

\begin{lemma}
Let $\varphi \in {\mathcal V}({\mathbb R}^n)$.  
Then the system $\{\exp \langle \lambda, z \rangle \}_{\lambda \in {\mathbb C}^n}$ is complete in 
$F^2_{\varphi}$.
\end{lemma}

{\bf Proof}. Let $S$ be a linear continuous functional on the space 
$F^2_{\varphi}$  such that  $S(e^{\langle \lambda , z \rangle}) = 0$ for each $\lambda \in {\mathbb C}^n$. Since for each multi-index 
$\alpha \in {\mathbb Z}_+^n$ we have that 
$
(D^{\alpha}_{\lambda} \hat S)(\lambda) = 
S(z^{\alpha} e^{\langle z, \lambda \rangle}),
$
then from this equality we get that
$S(z^{\alpha}) = 0$.
Since the function $\varphi(\vert z_1 \vert, \ldots , \vert z_n \vert)$ is convex on ${\mathbb C}^n$ then from the result by B.A. Taylor on the weight approximation of entire functions by polynomials \cite [Theorem 2] {T} it follows that the polynomials are dense in $F^2_{\varphi}$. Hence, $S$ is the zero functional. By the known corollary of Khan-Banach theorem we obtain that the system
$\{\exp \langle \lambda, z \rangle \}_{\lambda \in {\mathbb C}^n}$ is complete in $F^2_{\varphi}$. $\square$

Note that the system $\{z^{\alpha}\}_{\vert \alpha \vert \ge 0}$ is orthogonal in
$F^2_{\varphi}$. Besides that it is complete in $F^2_{\varphi}$. Therefore,
the system $\{z^{\alpha}\}_{\vert \alpha \vert \ge 0}$ is a basis in $F^2_{\varphi}$.

\begin{lemma}
Let $\varphi \in {\mathcal V}({\mathbb R}^n)$.  Then 
$$
c_{\alpha}(\varphi) \ge  \frac {\pi^n} {\tilde \alpha_1 \cdots \tilde \alpha_n}
e^{2 (\varphi [e])^*(\tilde \alpha)}, \ \alpha \in {\mathbb Z}_+^n.
$$
In particular, for each $M > 0$  there exists a constant $C_M > 0$ such that 
$c_{\alpha}(\varphi) \ge  C_M M^{\vert \alpha \vert}$
for each $\alpha \in {\mathbb Z}_+^n$.
\end{lemma}  

{\bf Proof}.
For each $\alpha \in {\mathbb Z}_+^n$ and each positive numbers $R_1, \ldots , R_n$ we have that
$$
c_{\alpha}(\varphi) = (2 \pi)^n \int \limits_{0}^{\infty} \cdots \int \limits_{0}^{\infty} r_1^{2 \alpha_1 + 1} \cdots 
r_n^{2 \alpha_n + 1} e^{-2\varphi(r_1, \cdots , r_n)} \, d r_1 \cdots d r_n \ge  
$$
$$
\ge (2 \pi)^n \int \limits_{0}^{R_1} \cdots \int \limits_{0}^{R_n} r_1^{2 \alpha_1 + 1} \cdots 
r_n^{2 \alpha_n + 1} e^{-2\varphi(R_1, \cdots , R_n)} \, d r_1 \cdots d r_n = 
$$
$$
= (2 \pi)^n \frac {R_1^{2 \alpha_1 + 2}} {2 \alpha_1 + 2} \cdots \frac {R_n^{2 \alpha_n + 2}} {2 \alpha_n + 2} e^{-2\varphi(R_1, \cdots , R_n)}.
$$
This implies that for each $t \in {\mathbb R}^n$ 
$$
c_{\alpha}(\varphi) \ge  \frac {\pi^n} {\tilde \alpha_1 \cdots \tilde \alpha_n}
e^{\langle 2 \tilde \alpha, t \rangle - 2\varphi [e](t)}.
$$
Hence,
$$
c_{\alpha}(\varphi) \ge  \frac {\pi^n} {\tilde \alpha_1 \cdots \tilde \alpha_n}
e^{2 (\varphi [e])^*(\tilde \alpha)}.
$$
Now using Proposition 2 we get the second statement of the lemma. $\square$

\begin{lemma}
Assume that an entire in ${\mathbb C}^n$ function  
$f(z) = \sum \limits_{\vert \alpha \vert \ge 0} a_{\alpha} z^{\alpha}$ is in $F^2_{\varphi}$. Then 
$
\sum\limits_{\vert \alpha \vert \ge 0} \vert a_{\alpha} \vert^2 c_{\alpha}(\varphi) < \infty 
$                                                
and
$
\Vert f \Vert_{\varphi}^2 = \sum_{\vert \alpha \vert \ge 0} \vert a_{\alpha} \vert^2 c_{\alpha}(\varphi).
$

Vice versa, let the sequence $(a_{\alpha})_{\vert \alpha \vert \ge 0}$ of complex numbers $a_{\alpha}$ be such that the series
$\sum \limits_{\vert \alpha \vert \ge 0} \vert a_{\alpha} \vert^2 c_{\alpha}(\varphi)$ converges. Then
$f(z) = \sum \limits_{\vert \alpha \vert \ge 0} a_{\alpha} z^{\alpha} \in H({\mathbb C}^n)$.
Moreover, $f \in F^2_{\varphi}$.

\end{lemma}                               

{\bf Proof}. Let $f(z) = \sum \limits_{\vert \alpha \vert \ge 0} a_{\alpha} z^{\alpha}$ be an entire function in ${\mathbb C}^n$ from $F^2_{\varphi}$. Then
$$
\Vert f \Vert_{\varphi}^2 = \int_{{\mathbb C}^n} \vert f(z)\vert^2 
e^{- 2 \varphi (abs \ z)} \ d \lambda (z) = 
$$
$$
=
\int_{{\mathbb C}^n} \sum_{\vert \alpha \vert \ge 0} a_{\alpha} z^{\alpha} 
\sum_{\vert \beta \vert \ge 0} {\overline a_{\beta}} {\overline z}^{\beta} 
e^{- 2 \varphi (abs \ z)} \ d \mu_n (z) = 
$$
$$
=
\sum_{\vert \alpha \vert \ge 0} \vert a_{\alpha} \vert^2 
\int_{{\mathbb C}^n} \vert z_1 \vert^{2 \alpha_1} \cdots \vert z_n \vert^{2 \alpha_n} 
e^{- 2 \varphi (abs \ z)} \ d \mu_n (z) = 
\sum_{\vert \alpha \vert \ge 0} \vert a_{\alpha} \vert^2 c_{\alpha}(\varphi). 
$$

Vice versa, the convergence of the series $\sum \limits_{\vert \alpha \vert \ge 0} \vert a_{\alpha} \vert^2 c_{\alpha}(\varphi)$ and Lemma 2 implies that for each $\varepsilon > 0$ there exists a constant $c_{\varepsilon} > 0$ such that  
$\vert a_{\alpha} \vert \le c_{\varepsilon} \varepsilon^{\vert \alpha \vert}$ 
for each $\alpha \in {\mathbb Z}_+^n$.
This means that $f(z) = \sum \limits_{\vert \alpha \vert \ge 0} a_{\alpha} z^{\alpha}$ is an entire function on ${\mathbb C}^n$. It is easy to see that $f \in F^2_{\varphi}$. \ $\square$

\begin{lemma}
Let $\varphi \in {\mathcal V}({\mathbb R}^n)$.  Then
$$
\frac {(2 \pi)^n} {e} 
V(D_{\tilde \alpha}^{\varphi [e]} (1 / 2)) e^{2(\varphi [e])^*(\tilde \alpha)} \le 
c_{\alpha}(\varphi) \le (2 \pi)^n (1 + n!) V(D_{\tilde \alpha}^{\varphi [e]} (1 / 2)) e^{2(\varphi [e])^*(\tilde \alpha)}
$$
for each $\alpha \in {\mathbb Z}_+^n$. 
\end{lemma}  

{\bf Proof}. Let  
$\alpha = (\alpha_1, \ldots , \alpha_n) \in {\mathbb Z}_+^n$. Then
$$
c_{\alpha}(\varphi) = (2 \pi)^n \int \limits_{0}^{\infty} \cdots \int \limits_{0}^{\infty} r_1^{2 \alpha_1 + 1} \cdots 
r_n^{2 \alpha_n + 1} e^{-2\varphi(r_1, \cdots , r_n)} \, d r_1 \cdots d r_n = 
$$
$$
= (2 \pi)^n \int \limits_{-\infty}^{\infty} \cdots \int \limits_{-\infty}^{\infty} e^{(2 \alpha_1 + 2) t_1 + \cdots + 
(2\alpha_n + 2) t_n - 2\varphi [e](t_1, \ldots , t_n)} \, d t_1 \cdots d t_n .
$$
That is,
$$
c_{\alpha}(\varphi) =(2 \pi)^n \int_{{\mathbb R}^n} 
e^{\langle 2 \tilde \alpha, t \rangle - 2\varphi [e](t)} \, d t.
$$
By Theorem B we have that 
$$
(2 \pi)^n e^{-1} V(D_{2 \tilde \alpha}^{2 \varphi [e]}) e^{2(\varphi [e])^*(\tilde \alpha)} \le 
c_{\alpha}(\varphi) \le (2 \pi)^n (1 + n!) V(D_{2 \tilde \alpha}^{2 \varphi [e]}) e^{2(\varphi [e])^*(\tilde \alpha)}
$$
Since $D_{2 \tilde \alpha}^{2 \varphi [e]} = D_{\tilde \alpha}^{\varphi [e]} \left(\frac 1 2\right)$ then from the last inequality our assertion follows. $\square$

{\bf 3.2. Proof of the Theorem}. 
Let us show that the mapping ${\mathcal L}$ acts from $(F^2_{\varphi})^*$ into $F^2_{\varphi^*}$.
Let $S \in (F^2_{\varphi})^*$. There exists a function 
$g_S \in F^2_{\varphi}$ such that 
$S(f) = (f, g_S)_{\varphi}$, that is, 
$$
S(f) = \int_{{\mathbb C}^n} f(z) \overline {g_S(z)} 
e^{- 2 \varphi (abs \ z)} \ d \mu_n (z), \ f \in F^2_{\varphi}.
$$
At that, $\Vert S \Vert = \Vert g_S \Vert_{\varphi}$.
If 
$
g_S(z) = \sum \limits_{\vert \alpha \vert \ge 0} b_{\alpha} z^{\alpha}
$
then
$
\hat S (\lambda) = \sum_{\vert \alpha \vert \ge 0} 
\frac {c_{\alpha}(\varphi) \overline {b_{\alpha}}} {\alpha!} \lambda^{\alpha}, \ \lambda \in  {\mathbb C}^n.
$
Hence, 
\begin{equation}
\Vert \hat S \Vert_{\varphi^*}^2 = 
\sum_{\vert \alpha \vert \ge 0} \left(\frac {c_{\alpha}(\varphi) \vert b_{\alpha}\vert}{\alpha!}\right)^2 c_{\alpha}(\varphi^*).
\end{equation}
By Lemma 3 
$$
c_{\alpha}(\varphi) \le (2 \pi)^n (1 + n!) V(D_{\tilde \alpha}^{\varphi [e]}(1 / 2)) e^{2(\varphi [e])^*(\tilde \alpha)},
$$
$$
c_{\alpha}(\varphi^*) \le 
(2 \pi)^n (1 + n!) V(D_{\tilde \alpha}^{\varphi^* [e]}(1 / 2)) e^{2(\varphi^*[e])^*(\tilde \alpha)}
$$
for each $\alpha \in {\mathbb Z}_+^n$.
Therefore, 
$$
c_{\alpha}(\varphi) c_{\alpha}(\varphi^*) \le 
$$
$$
\le 
(2 \pi)^{2n} (1 + n!)^2 
V(D_{\tilde \alpha}^{\varphi [e]}(1 / 2)) V(D_{\tilde \alpha}^{\varphi^* [e]}(1 / 2))
e^{2((\varphi [e])^*(\tilde \alpha) + (\varphi^*[e])^*(\tilde \alpha))}
$$
for each $\alpha \in {\mathbb Z}_+^n$.
According to Proposition 6 for each  
$\alpha = (\alpha_1, \ldots , \alpha_n) \in {\mathbb Z}_+^n$ we have that
$$
(\varphi [e])^*(\tilde \alpha) + (\varphi^* [e])^*(\tilde \alpha) = \sum \limits_{j =1}^n ((\alpha_j+1) \ln (\alpha_j+1) - (\alpha_j+1)).
$$
Since by the Stirling's formula \cite[P. 792] {Fiht} for each $m \in {\mathbb Z}_+$ we have that
$$
(m + 1) \ln (m+1) - (m + 1) = \ln \varGamma (m+1) - \ln \sqrt {2 \pi} + \frac 1 2 \ln (m+1) - \frac {\theta}{12(m+1)} , 
$$ 
where $\theta \in (0, 1)$ depends on $m$, then
$$
(\varphi [e])^*(\tilde \alpha) + (\varphi^* [e])^*(\tilde \alpha) = 
$$
$$
=-n  \ln \sqrt {2 \pi} + \sum \limits_{j =1}^n (\ln \varGamma (\alpha_j+1) + \frac 1 2 \ln (\alpha_j+1) - \frac {\theta_j}{12 (\alpha_j+1)}),
$$
where $\theta_j \in (0, 1)$ depends on $\alpha_j$.
Then
\begin{equation}
\frac {e^{2((\varphi [e])^*(\tilde \alpha) + (\varphi^* [e])^*(\tilde \alpha))}}{\alpha!^2} = \frac {1}{(2 \pi)^n} 
\prod \limits_{j =1}^n (\alpha_j + 1) e^{- \frac {\theta_j}{6 (\alpha_j+1)}} .
\end{equation}
Thus, 
$$
\frac {c_{\alpha}(\varphi) c_{\alpha}(\varphi^*)}{\alpha!^2} \le (2 \pi)^n (1 + n!)^2 
V(D_{\tilde \alpha}^{\varphi [e]}(1 / 2)) V(D_{\tilde \alpha}^{\varphi^* [e]}(1 / 2))
\prod \limits_{j =1}^n \tilde \alpha_j.
$$
Using the condition on $\varphi$ we get that  
$$
\frac {c_{\alpha}(\varphi) c_{\alpha}(\varphi^*)}{\alpha!^2} \le (2 \pi)^n (1 + n!)^2 K
$$
for each $\alpha  \in {\mathbb Z}_+^n$.
From this and (1) we obtain that
$$
\Vert \hat S \Vert_{\varphi^*}^2 \le M_1
\sum_{\vert \alpha \vert \ge 0} 
c_{\alpha}(\varphi) \vert b_{\alpha}\vert^2 = 
M_1 \Vert g_S \Vert_{\varphi}^2 = M_1 \Vert S \Vert^2,
$$
where $M_1 = (2 \pi)^n (1 + n!)^2 K$. Hence,  $\hat S \in F^2_{\varphi^*}$. Moreover, the latter estimate implies that the linear mapping ${\mathcal L}$ acts continuously from $(F^2_{\varphi})^*$ into $F^2_{\varphi^*}$. 

Note that the mapping ${\mathcal L}$ is injective since by Lemma~1 the system $\{\exp \langle \lambda, z \rangle \}_{\lambda \in {\mathbb C}^n}$ is complete in $F^2_{\varphi}$.

Let us show that the mapping ${\mathcal L}$ acts from $(F^2_{\varphi})^*$ onto $F^2_{\varphi^*}$. Let $G \in F^2_{\varphi^*}$. 
Using the representation of an entire function $G$ by the Taylor series 
$
G(\lambda) = \sum \limits_{\vert \alpha \vert \ge 0} d_{\alpha} \lambda^{\alpha}, \ \lambda \in {\mathbb C}^n,
$
we have
$
\Vert G \Vert_{\varphi^*}^2 = 
\sum \limits_{\vert \alpha \vert \ge 0} \vert d_{\alpha} \vert^2 c_{\alpha}(\varphi^*).
$
For each $\alpha  \in {\mathbb Z}_+^n$ define the numbers
$
g_{\alpha} = \frac {\overline {d_{\alpha}}\alpha!}{c_{\alpha}(\varphi)}
$
and consider the convergence of the series 
$\sum \limits_{\vert \alpha \vert \ge 0} \vert g_{\alpha} \vert^2 c_{\alpha}(\varphi)$.
We have that
$$
\sum \limits_{\vert \alpha \vert \ge 0} 
\vert g_{\alpha} \vert^2 c_{\alpha}(\varphi) = 
\sum \limits_{\vert \alpha \vert \ge 0} 
\left \vert \frac {\overline {d_{\alpha}}\alpha!}{c_{\alpha}(\varphi)} \right\vert^2 c_{\alpha}(\varphi) = \sum \limits_{\vert \alpha \vert \ge 0} 
\frac {\alpha!^2}{c_{\alpha}(\varphi)c_{\alpha}(\varphi^*)} 
\vert d_{\alpha}\vert^2  c_{\alpha}(\varphi^*).
$$
By Lemma 4 
$$
c_{\alpha}(\varphi) \ge e^{-1} V(D_{\tilde \alpha}^{\varphi [e]}(1 / 2)) e^{2(\varphi [e])^*(\tilde \alpha)},
$$
$$
c_{\alpha}(\varphi^*) \ge 
e^{-1} V(D_{\tilde \alpha}^{\varphi^* [e]}(1 / 2)) e^{2(\varphi^*[e])^*(\tilde \alpha)}
$$
for each $\alpha  \in {\mathbb Z}_+^n$.
Therefore,
$$
c_{\alpha}(\varphi) c_{\alpha}(\varphi^*) \ge e^{-2} V(D_{\tilde \alpha}^{\varphi [e]}(1 / 2)) V(D_{\tilde \alpha}^{\varphi^* [e]}(1 / 2)) e^{2((\varphi [e])^*(\tilde \alpha) + (\varphi^*[e])^*(\tilde \alpha))}
$$
for each $\alpha  \in {\mathbb Z}_+^n$.
From this and the equality (2) we have that
$$
\frac {\alpha!^2}{c_{\alpha}(\varphi) c_{\alpha}(\varphi^*)} \le  
\frac {e^2 (2 e \pi)^n}{V(D_{\tilde \alpha}^{\varphi [e]}(1 / 2)) V(D_{\tilde \alpha}^{\varphi^* [e]}(1 / 2)) 
\prod \limits_{j =1}^n (\alpha_j + 1)} 
$$
for each 
$\alpha = (\alpha_1, \ldots , \alpha_n) \in {\mathbb Z}_+^n$. 
Using the condition on $\varphi$ we obtain that 
$$
\frac {\alpha!^2}{c_{\alpha}(\varphi) c_{\alpha}(\varphi^*)} \le  
K e^2 (2 e \pi)^n, \ \forall \alpha  \in {\mathbb Z}_+^n.
$$
Therefore, for the considered series we have that
\begin{equation}
\sum \limits_{\vert \alpha \vert \ge 0} 
\vert g_{\alpha} \vert^2 c_{\alpha}(\varphi) \le K e^2 (2 e \pi)^n 
\sum \limits_{\vert \alpha \vert \ge 0} 
\vert d_{\alpha}\vert^2  c_{\alpha}(\varphi^*) = K e^2 (2 e \pi)^n \Vert G \Vert_{\varphi^*}^2.
\end{equation}
Thus, the series 
$\sum \limits_{\vert \alpha \vert \ge 0} \vert g_{\alpha} \vert^2 c_{\alpha}(\varphi)$ converges. But then by Lemma 3 the function
$$
g(\lambda) = \sum \limits_{\vert \alpha \vert \ge 0} g_{\alpha} \lambda^{\alpha}, \ \lambda \in {\mathbb C}^n,
$$
is entire and by (3) $g$ belongs to $F^2_{\varphi}$ and
\begin{equation}
\Vert g \Vert_{\varphi}^2 \le K e^2 (2 e \pi)^n \Vert G \Vert_{\varphi^*}^2.
\end{equation}
Define a functional $S$ on $F^2_{\varphi}$ by the formula
$$
S(f) = \int_{{\mathbb C}^n} f(z) \overline {g(z)} 
e^{- 2 \varphi (abs z)} \ d \mu_n (z), \ f \in F^2_{\varphi}.
$$
Obviously, $S$ is a linear continuous functional on $F^2_{\varphi}$ and
$\hat S = G$. Since $\Vert S \Vert = \Vert g \Vert_{\varphi}$ then estimate (4) shows that the inverse mapping ${\mathcal L}^{-1}$ is continuous. Thus,  ${\mathcal L}$ is an isomorphism between the spaces $(F^2_{\varphi})^*$ and $F^2_{\varphi^*}$. $\square$

Musin I. Kh.

Institute of Mathematics, Ufa Scientific Center, RAS,

Chernyshevsky str. 112,

450008, Ufa, Russia

E-mail: musin\_ildar@mail.ru


\bigskip
\begin{thebibliography}{1}

\bibitem {M} 
Il'dar Kh. Musin, {\it On a space of entire functions rapidly decreasing on $R^n$ and its Fourier transform},  Concrete Operators, {\bf 2}:1 (2015), 120-138

\bibitem {A 1} Daniel Azagra, {\it Global and fine approximation of convex functions}, Proc. London Math. Soc., {\bf 107}:4 (2013), 799-824

\bibitem {A 2}  D. Azagra, {\it Global approximation of convex functions},  arXiv:1112.1042v7.

\bibitem {T} B.A. Taylor, {\it On weighted polynomial approximation of entire functions}, Pacific Journal of Mathematics, {\bf 36}:2 (1971), 523-539

\bibitem {N-P} V.V. Napalkov, S.V. Popenov, {\it On the Laplace transform of functionals in the Bergman weight space of entire functions in ${\mathbb C}^n$},  Dokl. Akad. Nauk., {\bf 352}:5 (1997), 595-597; [Dokl. Math., {\bf 55}:1 (1997), 110-112]

\bibitem {P} S.V. Popenov, {\it On Laplace transform of functionals in some weighted Bergman spaces in ${\mathbb C}^n$}, Proc. Int. Conf. "`Complex analysis, differential equations, numerical methods and applications"', Ufa, {\bf 2} (1996), 125-132. (in Russian)

\bibitem {Yulm} R.S. Yulmukhametov, {\it Asymptotics of multidimensional Laplace integral}, "`Studies in approximation theory"', Ufa, 132-137 (1989). (in Russian)

\bibitem {NBY} V.V. Napalkov, R.A. Bashmakov, R.S.Yulmukhametov, {\it Asymptotic behavior of Laplace integrals and geometric characteristics of convex functions}, Dokl. Akad. Nauk., {\bf 413}:1 2007), 20-22; [Dokl.Math., {\bf 75}:2 (2007), 190-192]

\bibitem {BIY}
R.A. Bashmakov, K.P. Isaev, R.S. Yulmukhametov, {\it On geometric characteristics of convex function and Laplace integrals}, Ufimskij Matem. Zhurn., {\bf 2}:1 (2010), 3-16. (in Russian)

\bibitem {Fiht}
G.M. Fikhtengolts, {\it Course of differential and integral calculus}, V. II, Nauka, Moscow (1970). (in Russian)

\end{thebibliography}
\end{document}